\documentclass[3p,preprint]{elsarticle}

\usepackage{amsmath}
\usepackage{amsfonts}
\usepackage{amssymb,latexsym}
\usepackage{amsthm}
\usepackage{graphicx}

\newtheorem{theorem}{Theorem}[section]

\newtheorem{Remark}{Remark}[section]
\newtheorem{Corollary}{Corollary}[section]

\usepackage{color}
\usepackage{newlfont}
\usepackage{subfigure}
\usepackage{verbatim}
\usepackage{rotating}
\usepackage{multirow}
\usepackage{units}
\usepackage{bm}
\usepackage[english]{babel}
\usepackage[utf8]{inputenc}
\usepackage{algorithm}
\usepackage[noend]{algpseudocode}
\usepackage{booktabs}
\usepackage{caption}
\usepackage{hyperref}

\newcommand{\red}[1]{{\color{black} #1}}

\journal{}

\begin{document}

\begin{frontmatter}

\title{Random sampling and unisolvent interpolation \\ by almost everywhere analytic functions}



\author[address-CA,address-GNCS]{Francesco Dell'Accio}
\ead{francesco.dellaccio@unical.it}

\author[address-PD,address-GNCS]{Alvise Sommariva}
\ead{alvise@math.unipd.it}

\author[address-PD,address-GNCS]{Marco Vianello\corref{corrauthor}}
\ead{marcov@math.unipd.it}

\cortext[corrauthor]{Corresponding author}

\address[address-CA]{University of Calabria, Cosenza, Italy}

\address[address-PD]{University of Padova, Italy}

\address[address-GNCS]{Member of the INdAM Research group GNCS}

\begin{abstract}
We prove {\em a.s.} (almost sure) unisolvency of interpolation by continuous random sampling with respect to any given density, in spaces of multivariate {\em a.e.} (almost everywhere) analytic functions. Examples are given concerning polynomial and RBF approximation. 
\end{abstract}



\end{frontmatter}

\vskip1cm

This note is aimed at proving a general result on unisolvency of multivariate interpolation in spaces of {\em a.e.} (almost everywhere) analytic functions, by continuous random sampling with any density. 
Such a result might be already known in some specific function space, and we do not even attempt to give a partial overview of the vast literature on random sampling and its connections with approximation theory. On the contrary, we wish to emphasize that in the random setting there is a common framework for apparently quite far approaches, such as for example polynomial approximation and RBF approximation. 

Such a common framework is given by real analytic functions. We recall that a  function is real analytic on an open subset $\Omega\subset \mathbb{R}^d$ if for each $x\in \Omega$ the function may be represented by a convergent power series in some neighborhood of $x$; cf. \cite[Ch.2]{KP02}.

\begin{theorem}
Let $\Omega$ be an open connected subset of $\mathbb{R}^d$ and $\{f_j\}_{j\geq 1}$ a set of functions defined on $\Omega$, such that 
\begin{itemize}
\item[$i)$] each $f_j$ is real analytic up to a set of null Lebesgue measure, say $I_j\subset \Omega$;
\item[$ii)$] the $\{f_j\}_{j\geq 1}$ are linearly independent on every connected component of $\Omega\setminus I$, where $I=\bigcup_{j=1}^m{I_j}$. 
\end{itemize}
\noindent
Moreover, let $\{x_i\}_{i\geq 1}$ be a randomly distributed sequence on $\Omega$ with respect to any given probability density $\sigma(x)$, i.e. a point sequence produced by sampling a sequence of continuous random variables $\{X_i\}_{i\geq 1}$ which are independent and identically distributed in $\Omega$ with density $\sigma(x)$.  
Then, for every $m\geq 1$ the matrix $V_m=[f_j(x_i)]$, $1\leq i,j\leq m$, is {\em a.s.} (almost surely) nonsingular.

\end{theorem}
\vskip0.2cm
\noindent{\bf Proof.} We proceed by induction on $m$. 

First, we prove the assertion for $m=1$. In fact, $det(V_1)=f_1(x_1)=0$ iff $x_1$ falls on the zero set of 
 $det(f_1(x))=f_1(x)$, $x\in \Omega$. Now, the zero set of $f_1$ in $\Omega$, say $\mathcal{Z}(f_1)$, is the disjoint union of the zero set in $\Omega\setminus I$, say $\mathcal{Z}_1(f_1)$, with the possible zero set in $I$, say $\mathcal{Z}_2(f_1)$. Observe that $f_1$ is real analytic on the set $\Omega\setminus I$, which is open since $I$ is closed 
 in the induced topology of $\Omega$, each $f_j$ being analytic in $\Omega \setminus I_j$. 
Moreover, by assumption $f_1$ is not identically zero on each connected component of $\Omega\setminus I$. Notice that such connected components are at most countable. 

 Hence the zero set $\mathcal{Z}_1(f_1)$, which is the union of the zero sets of the connected components of $\Omega\setminus I$, has null Lebesgue measure by a well-known basic result of measure theory, asserting that the zero set of a nonzero real analytic function on an open connected set in $\mathbb{R}^d$ has null Lebesgue measure (cf. \cite{M20} for an elementary proof). 
 
 In turn, the zero set $\mathcal{Z}_2(f_1)$  has null Lebesgue measure being a subset of $I$, and thus $\mathcal{Z}(f_1)$ has null Lebesgue measure. Consequently, $\mathcal{Z}(f_1)$ has null measure also with respect to the density $\sigma(x)$, since $\int_{\mathcal{Z}(f_1)}{\sigma(x)\,dx}=0$, i.e. $f_1(x_1)$ 
 is {\em a.s.} nonzero (i.e., with probability 1). 

Now, let assume that $V_m$ is nonsingular with probability 1, and consider the $(m+1) \times (m+1)$ matrix $U_{m+1}(x)$ obtained by adding to $V_m$ the $(m+1)$-th column $$[f_{m+1}(x_1),\dots,f_{m+1}(x_m),f_{m+1}(x)]^t$$ and the $(m+1)$-th row 
$$[f_1(x),\dots,f_m(x),f_{m+1}(x)]\;.$$
Applying Laplace rule to the last row, we get that 
$$
det(U_{m+1}(x))=det(V_m)f_{m+1}(x)+\alpha_mf_m(x)+\dots
+\alpha_1f_1(x)
$$
where $\alpha_1,\dots,\alpha_m$ are the other corresponding minors with the 
appropriate sign. Notice that $det(U_{m+1}(x))$ is not identically zero on each connected component of $\Omega\setminus I$ with probability 1, since the $\{f_j\}$ are linearly independent there, and $det(V_m)\neq 0$ with probability 1, by inductive hypothesis. 

Thus, by the same arguments of the $m=1$ instance with $x_{m+1}$ and $det(U_{m+1}(x))$ substituting $x_1$ and $f_1(x)$, respectively, we get that 
the probability that $x_{m+1}$ falls on the zero set in $\Omega$ 
of $det(U_{m+1}(x))$ is null, i.e. $V_{m+1}=U_{m+1}(x_{m+1})$ is nonsingular with probability 1. To be more precise, $$prob\{det(V_{m+1})=0\}= prob\{det(V_{m+1})=0\;\&\; det(V_m)=0\}$$
$$+\,prob\{det(V_{m+1})=0\;\&\;det(V_m)\neq 0\}=0+0=0\;,$$
since the events are disjoint and both their probabilities are null (notice that the first one has null probability since it is a subevent of the event $det(V_m)=0$, which has null probability by inductive hypothesis).
\hspace{0.2cm} $\square$
\vspace{0.5cm}

\begin{Remark} 
{\em 
We stress that the above result is valid with {\em any distribution density}. For example, in a box, points can have a uniform distribution, but also a normal or a product Chebyshev distribution. It is also worth noticing that the case of everywhere analytic functions has been recently considered in \cite{XN23}, within an abstract context named ``$\mu ZC$ sequences'' of finite-dimensional function spaces (see also the references therein for previous univariate results).
}
\end{Remark}

\begin{Remark}  
{\em 
 Theorem 0.1 concerns in principle multivariate interpolation by random sampling. 
 On the other hand, {\em a.s.} (almost sure) unisolvency at random nodes is relevant also in the least-squares framework, since the corresponding rectangular Vandermonde-like matrix $V$ has {\em a.s.} maximal rank, because it contains an {\em a.s.} nonsingular interpolation matrix and the Gram matrix $V^tV$ becomes {\em a.s.} definite positive.
 }
\end{Remark}

Below, to the purpose of illustration, we give some applications of the above result in different function spaces.
As a first application, we observe that {\em point sequences of the appropriate length, randomly distributed with respect to any given probability density on an open connected set, are {\em a.s.} unisolvent 
for multivariate interpolation by polynomials, trigonometric polynomials, and rational functions whose denominator does not vanish on the domain}. 
Indeed, it is sufficient to observe that both polynomials and trigonometric polynomials are entire functions, whereas rational functions are real analytic. 
To quote some examples among many others in polynomial approximation, unisolvency at random nodes is useful within different topics such as numerical differentiation by local polynomial interpolation, multinode Shepard-like interpolation, compressed (Quasi)MonteCarlo integration, randomized weakly admissible meshes and polynomial least-squares; cf. e.g. \cite{DADT19,DADTSV22,ESV22,XN23} with the references therein.
\vskip0.2cm

We can now give a relevant Corollary of Theorem 0.1, concerning RBF interpolation.

\begin{Corollary}
Sequences $x_1,\dots,x_m$ of randomly distributed points with respect to 
any given probability density on an open connected set $\Omega\subset \mathbb{R}^d$, are {\em a.s.} unisolvent 
for multivariate RBF interpolation with fixed 
distinct centers $\{\xi_1,\dots,\xi_m\}\subset \Omega$ by Gaussians, Multiquadrics (MQ), Inverse Multiquadrics (IMQ), and Thin-Plate Splines (TPS) for $d
\geq 1$, and also by Radial Powers (RP) for $d\geq 2$.
\end{Corollary}
\vskip0.2cm
\noindent{\bf Proof.}
The quoted RBF are of the form $f_j(x)=\phi(\|x-\xi_j\|_2)$ with univariate radial functions, respectively, $\phi(r)=e^{-r^2}$ (Gaussians), $\phi(r)=(1+r^2)^{1/2}$ (MQ), $\phi(r)=(1+r^2)^{-1/2}$ (IMQ), $\phi(r)=r^k$, $k$ odd (RP), and $\phi(r)=r^k\log(r)$, $k$ even (TPS); cf. e.g. \cite{F07}.  

Gaussians, MQ and IMQ are real analytic in $\mathbb{R}^d$, since the corresponding $\phi(r)$ is real analytic in $\mathbb{R}$ and the composition of real analytic functions is real analytic \cite[Prop.1.4.2, p.19]{KP02}. {Moreover, they are linearly independent on $\Omega$. Indeed, Gaussians and IMQ are strictly positive definite and the corresponding standard interpolation matrix at the centers, $U=[\phi(\|\xi_i-\xi_j\|_2]$, $1\leq i,j\leq m$, is positive definite; cf. e.g. \cite{F07}. On the other hand, MQ are linearly independent in view of a classical result of Micchelli on the invertibility of $U$ 
for conditionally positive definite RBF of order 1; cf. \cite{M86}}. Then Theorem 0.1 applies. 

On the other hand, RP and TPS are everywhere continuous, and  
real analytic with the exception of the center $\xi_j$, since in both cases $\phi(\sqrt{\cdot}\,)$ is real analytic in $\mathbb{R}^+$. This fact also ensures that for $d\geq 2$ they are linearly independent on the unique connected component of $\Omega\setminus \{\xi_1,\dots, \xi_m\}$. Indeed, if they were dependent there, by continuity they would be dependent also on the whole $\Omega$. But this not possible (for any $d\geq 1$ in fact), since 
if they were dependent, one of them would be a linear combination of the others, and would then result analytic at its center. Again, Theorem 0.1 applies. 

In the case of univariate TPS, $f_j(x)=|x-\xi_j|^k\log(|x-\xi_j|)$, $\Omega=(a,b)$ 
and $a<\xi_1<\xi_2<\dots<\xi_m<b$. Assume that a linear combination 
$f(x)=\alpha_1 f_1(x)+\dots +\alpha_m f_m(x)\equiv 0$ in $(\xi_\ell,\xi_{\ell+1})$, 
where we set $\xi_{-1}=a$ and $\xi_{m+1}=b$. Now, since the $(k+1)$-th derivatives $f_\ell^{(k+1)}(x)$ and $f_{\ell+1}^{(k+1)}(x)$ tend to $\infty$ as $x\to \xi_\ell^+$ and $x\to \xi_{\ell+1}^-$ respectively, due to the presence of the logarithmic factor, necessarily $\alpha_\ell=0=\alpha_{\ell+1}$, otherwise   
$f^{(k+1)}(x)$ would tend to $\infty$ at $\xi_\ell$ and $\xi_{\ell+1}$. Then, $f(x)$ being analytic in 
$(\xi_{\ell-1},\xi_{\ell+2})$ and identically zero in the subinterval $(\xi_\ell,\xi_{\ell+1})$, it is $f(x)\equiv 0$ in the whole interval $(\xi_{\ell-1},\xi_{\ell+2})$ (where we set for convenience $\xi_j=a$ 
for $j<1$ and $\xi_j=b$ for $j>m$).
Repeating the reasoning above on all the progressively enlarging subintervals, we then obtain that $\alpha_j=0$ for all $j$, 
and thus $f_1,\dots,f_m$ are linearly independent on every open subinterval corresponding to consecutive centers, that is on every connected component 
of $(a,b)\setminus \{\xi_1,\dots, \xi_m\}$, and Theorem 0.1 then applies.
\hspace{0.2cm} $\square$
\vskip0.2cm

\begin{Remark} 
{\em
In the univariate case, for $m>k+1$ univariate RP are certainly linearly dependent in every subinterval determined by consecutive centers, being polynomials of degree $k$ there, and consequently Theorem 0.1 does not apply.
}
\end{Remark}

\begin{Remark} 
{\em The fact that least-squares approximation by Thin-Plate Splines (without any polynomial augmentation) can be desiderable in certain applications, has been
recognized in the literature (cf. e.g. \cite{P22}). In this framework, it is important to have a nonsingular Gram matrix, see Remark 0.2. For the deterministic theory of least-squares RBF approximation we may quote e.g. \cite[Ch.20]{F07}, 
\cite[\S 3.10]{I04} and the classical papers \cite{QSW93,SW93}.

On the other hand, Corollary 0.1 suggests the following operative procedure for {\em interpolation by TPS and RP without any polynomial augmentation}. Given a sample at random points 
$x_1,\dots,x_m$ distributed with respect to any given probability density, we can draw another random point distribution of the same length, say $\xi_1,\dots,\xi_m$ (for example using a uniform distribution or even the same density, the probability that this distributions have common points being null). Then, we can interpolate at $x_1,\dots,x_m$ with the 
RBF basis $\{\phi(\|x-\xi_j\|_2)\}$, $1\leq j\leq m$, that is with the interpolation matrix $[\phi(\|x_i-\xi_j\|_2)]$, $1\leq i,j\leq m$, which will be {\em a.s.} nonsingular.

In alternative, we can fix a priori a set of distinct centers $\xi_1,\dots,\xi_m$ with any general or application driven criterion, and then proceed by continuous random sampling of length $m$. Again, the interpolation matrix will be 
{\em a.s.} nonsingular.
}
\end{Remark}
    
As a last example, we give the following corollary concerning in particular polynomial interpolation on surfaces.

\begin{Corollary}
Let $\mathcal{S}\subset \mathbb{R}^3$ be a surface that admits an analytic parametrization $x=\psi(u,v)$ from a connected open set $D\subset \mathbb{R}^2$, i.e. $\psi=(\psi_1,\psi_2,\psi_3)$ where $\psi_i:D\to \mathbb{R}^3$ are analytic 
and $\psi(D)=\mathcal{S}$. Moreover, let $\{p_j\}_{j\geq 1}$ be a set of trivariate polynomials that are linearly independent on $\mathcal{S}$ and $\{(u_i,v_i)\}_{i\geq 1}$ a randomly distributed sequence on $D$ with respect to any given probability density. 
Then the points $\{x_i=\psi(u_i,v_i)\}_{1\leq i\leq m}$ are {\em a.s.} unisolvent 
for polynomial interpolation in $\mbox{\em span}(p_1,\dots,p_m)$.
\end{Corollary}
\vskip0.2cm

\noindent{\bf Proof.} The proof is immediate in view of Theorem 0.1, by observing that the functions 
$f_j(u,v)=p_j(\psi(u,v))$ are linearly independent and real analytic on 
$\Omega=D$.
\hspace{0.2cm} $\square$
\vskip0.2cm 

{
For the purpose of illustration, relevant examples are regions of sphere, torus and cylinder, where there are natural entire 
parametrizations of trigonometric or algebraic/trigonometric type (spherical, toroidal and cylindrical coordinates), as well as Cartesian graphs of bivariate analytic functions. Some applications of Corollary 0.2 can be found, e.g., in the paper \cite{ESV23} that extends the method in \cite{ESV22} to (Quasi)Monte-Carlo integration on surfaces with a regular analytic parametrization, where it is required that the points be distributed with respect to the surface measure density, that is $\sigma(u,v)=\|\partial_u\psi\times \partial_v\psi\|_2/area(S)$. Another application arises in the context of multinode Shepard interpolation on the sphere with enhanced polynomial reproduction, cf. \cite{DADT23}.
}
\vskip0.5cm 
\noindent
{\bf Acknowledgements.} The authors wish to thank Len Bos for his illuminating  
observations concerning RBF spaces, and Federico Piazzon for his suggestions on relevant references.

Work partially supported by the DOR funds of the University of Padova, and by the INdAM-GNCS 2022 Projects “Methods and software for multivariate integral models” and ``Computational Methods for Kernel based Approximation and its Applications''. 

This research has been accomplished within the RITA ``Research ITalian network on Approximation" and the SIMAI Activity Group ANA\&A, and the UMI Group TAA ``Approximation Theory and Applications" (F. Dell'Accio, A. Sommariva).

\end{document}